\documentclass[12pt]{article}

\usepackage{a4,graphicx,amsmath,amsfonts,amssymb,bbm}
\usepackage{url}

\newcommand{\real}{\mathbbm{R}}

\newcommand{\ltwo}{\mathcal{L}^2(\Pi,\rho)}

\newtheorem{theorem}{Theorem}

\newtheorem{lemma}{Lemma}

\begin{document}
\begin{center}
{\bf \Large Model order reduction for \\[0.5ex] 
random nonlinear dynamical systems and \\[0.5ex] 
low-dimensional representations for \\[1ex] 
their quantities of interest}

\vspace{5mm}

{\large Roland~Pulch}

{\small Institut f\"ur Mathematik und Informatik,
  Universit\"at Greifswald, \\
  Walther-Rathenau-Str.~47, 17489 Greifswald, Germany. \\
Email: {\tt roland.pulch@uni-greifswald.de}}

\end{center}

\bigskip\bigskip


\begin{center}
{Abstract}

\begin{tabular}{p{13cm}}
We examine nonlinear dynamical systems of ordinary differential equations 
or differential algebraic equations. 
In an uncertainty quantification, 
physical parameters are replaced by random variables.
The state or inner variables as well as a quantity of interest
are expanded into 
series with orthogonal basis functions like the polynomial chaos expansions, 
for example.
On the one hand, the stochastic Galerkin method yields a large 
coupled dynamical system. 
On the other hand, a stochastic collocation method, 
which uses a quadrature rule or a sampling scheme, 
can be written in the form of a large weakly coupled dynamical system. 
We apply projection-based methods of nonlinear model order reduction 
to the large systems. 
A reduced-order model implies a low-dimensional representation of 
the quantity of interest. 
We focus on model order reduction by proper orthogonal decomposition. 
The error of a best approximation located in a low-dimensional subspace 
is analysed. 
We illustrate results of numerical computations for test examples. 

\bigskip

\end{tabular}
\end{center}

\clearpage


\section{Introduction}
Mathematical modelling often generates
dynamical systems of ordinary differential equations (ODEs) or 
differential algebraic equations (DAEs) in science and engineering. 
A quantity of interest (QoI) is defined by the state variables or 
inner variables. 
Multiple sources of uncertainties may affect the included 
parameters like modelling errors and measurement errors, 
for example. 
Several techniques have been developed to quantify the effects of
those uncertainties in the model predictions,
see~\cite{klir-wierman,sullivan,xiu-book}.
A common approach is the substitution of uncertain parameters
by random variables.  
Numerical methods often apply orthogonal expansions with 
unknown coefficient functions and predetermined 
random-dependent basis functions. 
The stochastic Galerkin method
 (intrusive method) 
produces a larger coupled dynamical system. 
A stochastic collocation technique
 (non-intrusive method) 
using a quadrature rule or a 
sampling method can be written as a single large weakly coupled 
dynamical system, see~\cite{pulch17}. 

A low-dimensional approximation of the random QoI often exists, 
where only a few basis functions are required for a sufficiently 
accurate representation.
In a specific context, the representation can be interpreted as a 
sparse approximation. 
Such low-dimensional representations have been computed by 
different methods like
least angle regression~\cite{blatman}, 
sparse grid quadrature~\cite{conrad-marzouk}, 
compressed sensing~\cite{doostan} 
and $\ell_1$-minimisation~\cite{jakeman2015,jakeman2016}.
Stochastic reduced bases were investigated for 
random linear systems of algebraic equations in~\cite{nair,sachdeva}.

Model order reduction (MOR) becomes favourable in the stochastic Galerkin 
method and the stochastic collocation approach
due to the high dimensionality of the systems.
In~\cite{pulch-matcom18}, it was shown that an MOR of the Galerkin system 
implies a low-dimensional approximation of the QoI in the case 
of linear dynamical systems. 
In~\cite{pulch-scee2016}, this strategy was carried over to an MOR 
of the single auxiliary system of the stochastic collocation method 
in the linear case.
In this paper, we consider nonlinear dynamical systems with random parameters,
where the QoI still depends linearly on 
state variables or inner variables. 
Thus an MOR generates a low-dimensional representation of the 
random QoI again.

More precisely, MOR is the tool to identify an approximation of the QoI
with a low number of basis functions. 
We do not investigate a reduction of the number of random parameters,
cf.~\cite{pulch-maten-augustin}.

Several efficient MOR methods exist for linear dynamical systems, 
see~\cite{antoulas,benner-mehrmann,gugercin,schilders}. 
Error bounds are available by Hardy norms of transfer functions in 
the frequency domain. 
Yet MOR for nonlinear dynamical systems still represents a challenging task. 
Typically, projection-based MOR schemes are applied like 
the proper orthogonal decomposition (POD), 
see~\cite{antoulas,kunisch}, 
or the trajectory piecewise linear approach, 
see~\cite{mohaghegh,rewienski}. 
We focus on the POD method, 
which identifies a low-dimensional representation of the QoI. 
Furthermore, the representation can be improved by the computation of 
a best approximation in the low-dimensional subspace. 
We prove an error bound for this best approximation. 
In the first place, our aim is the identification of a sufficiently 
accurate approximation with
as few basis functions as possible.  
In the second place, the MOR methods should decrease the 
computational effort.

The paper is organised as follows. 
Section~\ref{sec:problem-def} incloses the problem setup. 
Numerical methods yield large dynamical systems formulated
in Section~\ref{sec:methods}. 
The main part is given by Section~\ref{sec:mor}, 
where we apply MOR and show error bounds. 
We demonstrate the results of numerical experiments for two 
illustrative examples in Section~\ref{sec:example}.


\section{Problem definition}
\label{sec:problem-def}
We describe the task consisting in the identification 
of low-dimensional representations for QoIs from 
random nonlinear dynamical systems.

\subsection{Deterministic model}
Let a nonlinear dynamical system be given in the form
\begin{equation} \label{dae}
\begin{array}{rcl}  
E(p) \dot{x}(t,p) & = & 
A(p) x(t,p) + F(x(t,p),p) + B(p) u(t) \\[1ex] 
y(t,p) & = & C(p) x(t,p) \\
\end{array} 
\end{equation}
with matrices and functions depending on 
physical parameters $p \in \Pi \subseteq \real^q$. 
The sizes of the matrices are $A,E \in \real^{n \times n}$, 
$B \in \real^{n \times n_{\rm in}}$, $C \in \real^{n_{\rm out} \times n}$. 
The system involves a nonlinear function 
$F : \real^n \times \Pi \rightarrow \real^n$. 
For non-singular matrices~$E$, the system consists of ODEs 
with the state variables 
$x : [t_0,t_{\rm end}] \times \Pi \rightarrow \real^n$. 
For singular matrices~$E$, a system of DAEs is given with 
the inner variables~$x$. 
We consider initial value problems 
\begin{equation} \label{ivp}
x(t_0,p) = x_0(p) 
\qquad \mbox{for} \;\; p \in \Pi 
\end{equation}
with a predetermined function $x_0 : \Pi \rightarrow \real^n$. 
In the case of DAEs, the initial values have to satisfy consistency 
conditions and typically depend on the physical parameters of the system.

An input $u : [t_0,t_{\rm end}] \rightarrow \real^{n_{\rm in}}$ is supplied 
to the system~(\ref{dae}). 
An output $y : [t_0,t_{\rm end}] \times \Pi \rightarrow \real^{n_{\rm out}}$ 
is defined by the state variables or inner variables and the 
matrix~$C$. 
Without loss of generality, we restrict the analysis to the 
case of single-input-single-output (SISO) with $n_{\rm in} = n_{\rm out} = 1$. 
On the one hand, the results also apply to dynamical systems 
with a general nonlinear dependence of the right-hand side on the input. 
Moreover, the theory is applicable to autonomous dynamical systems. 
On the other hand, the linear dependence of the output on the 
state variables or inner variables is essential in this paper.

\subsection{Stochastic model} 
We adopt a common approach in uncertainty quantification (UQ),
see~\cite{sullivan,xiu-book}, for example.
Assuming that the parameters~$p \in \Pi$ of the system~(\ref{dae}) 
are uncertain, they are replaced by independent random variables 
$p : \Omega \rightarrow \Pi$ on some probability space 
$(\Omega,\mathcal{A},P)$ with event space~$\Omega$, 
sigma-algebra~$\mathcal{A}$ and probability measure~$P$.
A joint probability density function $\rho : \Pi \rightarrow \real$ 
is available in the case of traditional probability distributions.
For a measurable function $f : \Pi \rightarrow \real$, the expected value 
reads as
\begin{equation} \label{expectedvalue} 
\mathbb{E} \left[ f \right] := 
\int_{\Pi} f (p) \rho(p) \; \mbox{d}p 
\end{equation}
provided that the integral is finite.
The Hilbert space
$$ \ltwo := \left\{ f : \Pi \rightarrow \real \; : \; 
f \; \mbox{measurable and} \;
\mathbb{E} \left[ f^2 \right] < \infty \right\} $$
is equipped with the inner product
\begin{equation} \label{innerproduct} 
< f , g > \; := \mathbb{E} \left[ f g \right] =  
\int_{\Pi} f (p) g (p) \rho(p) \; \mbox{d}p 
\qquad \mbox{for} \;\; f,g \in \ltwo . 
\end{equation}
The accompanying norm reads as
$$ \left\| f \right\|_{\ltwo} := \sqrt{ < f , f> } . $$
Concerning the dynamical system~(\ref{dae}), 
we assume that 
$x_1(t,\cdot),\ldots,x_n(t,\cdot), y(t,\cdot) \in \ltwo$
pointwise for $t \in [t_0,t_{\rm end}]$. 

Now let a complete orthonormal system 
$( \Phi_i )_{i \in \mathbb{N}} \subset \ltwo$ 
be given. 
It holds that $< \Phi_i , \Phi_j > = 0$ for $i \neq j$ 
and $< \Phi_i , \Phi_j > = 1$ for $i = j$. 
We assume that the first basis function is always the constant function 
$\Phi_1 \equiv 1$.
It follows that the expansions 
\begin{equation} \label{pce} 
x(t,p) = \sum_{i=1}^\infty v_i(t) \Phi_i(p) 
\qquad \mbox{and} \qquad 
y(t,p) = \sum_{i=1}^\infty w_i(t) \Phi_i(p) ,
\end{equation}
where the coefficient functions 
$v_i : [t_0,t_{\rm end}] \rightarrow \real^n$ and 
$w_i : [t_0,t_{\rm end}] \rightarrow \real$ are defined by
\begin{equation} \label{coefficients}
v_{i,j}(t) = \; < x_j(t,\cdot) , \Phi_i(\cdot) > 
\qquad \mbox{and} \qquad 
w_i(t) = \; < y(t,\cdot) , \Phi_i(\cdot) > , 
\end{equation}
converge in $\ltwo$ pointwise for $t$ 
and $j=1,\ldots,n$
 with $x = (x_1,\ldots,x_n)^\top$.  

\subsection{Low-dimensional orthogonal representations}

Our aim is to identify a low-dimensional approximation using just
a few basis functions.
Several methods to obtain sparse or low-dimensional
representations have already been designed and investigated, see
\cite{blatman,conrad-marzouk,constantine-etal,doostan,jakeman2015,jakeman2016}.

In practice, the series~(\ref{pce}) have to be truncated to 
finite approximations.
A truncation yields 
\begin{equation} \label{pce-appr} 
x^{(m)}(t,p) = \sum_{i=1}^m v_i(t) \Phi_i(p) 
\qquad \mbox{and} \qquad 
y^{(m)}(t,p) = \sum_{i=1}^m w_i(t) \Phi_i(p) 
\end{equation}
for an integer $m \ge 1$. 
We investigate the output of the random dynamical system~(\ref{dae}) 
as QoI.
The truncation error reads as
\begin{equation} \label{truncationerror}
\left\| y(t,\cdot) - y^{(m)}(t,\cdot) 
\right\|_{\ltwo} = 
\sqrt{\sum_{i = m+1}^{\infty} w_i(t)^2} 
\end{equation}
for each~$t$ due to the orthonormality of the basis functions. 

Often orthonormal polynomials are chosen as basis functions 
with respect to the theory of the polynomial chaos (PC) expansions,

see~\cite{cameron-martin,ernst,xiu-book}. 
Therein, the multivariate basis polynomials are just the 
products of the univariate orthonormal polynomials. 
Each traditional probability distribution implies its 
family of univariate orthonormal polynomials, 
see~\cite{xiu-karniadakis}.
To obtain an initial set of basis functions, 
often all polynomials up to a total degree~$d$ are included.
Hence the number of basis functions is
$$ m = \frac{(q+d)!}{q! d!} , $$
see~\cite[p.~65]{xiu-book}.
Even if the total degree is moderate (say $3 \le d \le 5$), 
the number of basis polynomials becomes large for high numbers~$q$ 
of random parameters.

A numerical method yields approximations $\hat{w}_1,\ldots,\hat{w}_m$
of the exact coefficient functions in the 
truncated expansion~(\ref{pce-appr}). 
The investigated QoI reads as
\begin{equation} \label{qoi-method} 
\hat{y}(t,p) = \sum_{i=1}^m \hat{w}_i(t) \Phi_i(p) . 
\end{equation}
Assuming a large number~$m$, 
our aim is to identify a sufficiently accurate 
low-dimensional representation of~(\ref{qoi-method})
\begin{equation} \label{qoi-new} 
\bar{y}(t,p) = \sum_{i=1}^r \bar{w}_i(t) \Psi_i(p) 
\end{equation}
with a new orthonormal basis $\{ \Psi_1 , \ldots , \Psi_r \}$ satisfying 
$$ {\rm span} \{ \Psi_1 , \ldots , \Psi_r \} \subset 
{\rm span} \{ \Phi_1 , \ldots , \Phi_m \} $$
and new coefficient functions $\bar{w}_1 , \ldots , \bar{w}_r$ 
for some $r \ll m$. 
A special case is given by the selection of a subset
$\{ \Psi_1 , \ldots , \Psi_r \} \subset \{ \Phi_1 , \ldots , \Phi_m \}$. 
In this situation, the approximation~(\ref{qoi-new}) can be 
interpreted as a sparse representation, see~\cite{blatman}. 

The error of the complete approach is estimated by
$$ \begin{array}{rcl}
\left\| y(t,\cdot) - \bar{y}(t,\cdot) \right\|_{\ltwo} 
& \le & 
\left\| y(t,\cdot) - y^{(m)}(t,\cdot) \right\|_{\ltwo} \\[1ex]
& & \mbox{} + 
\left\| y^{(m)}(t,\cdot) - \hat{y}(t,\cdot) \right\|_{\ltwo} \\[1ex]
& & \mbox{} + 
\left\| \hat{y}(t,\cdot) - \bar{y}(t,\cdot) \right\|_{\ltwo}  \\
\end{array} $$
for each~$t$.
The upper bound consists of three terms: 
the truncation error, the error of the numerical method and 
the additional error of the low-dimensional approximation. 
We assume that the first and second term are sufficiently small 
due to choosing a large initial set of basis functions and a 
sufficiently accurate numerical method. 
The third term is analysed in Section~\ref{sec:erroranalysis}.


\section{Numerical Techniques}
\label{sec:methods}
We apply two well-known classes of numerical methods for the computation of 
the unknown coefficient functions~(\ref{coefficients}) in the orthogonal 
expansions~(\ref{pce}) now.

\subsection{Stochastic Galerkin method}
The stochastic Galerkin technique (intrusive method) 
represents a general approach, 
which can be applied to all types of differential equations 
including random variables, see~\cite{sullivan,xiu-book}. 
Linear or nonlinear dynamical systems were considered 
in~\cite{augustin-rentrop,pulch11,pulch13,pulch14}, for example. 
In our problem, the Galerkin method changes the 
nonlinear dynamical system~(\ref{dae}) into the larger coupled system
\begin{equation} \label{galerkin}
\begin{array}{rcl}
\hat{E} \dot{\hat{v}}(t) & = & 
\hat{A} \hat{v}(t) + \hat{F}(\hat{v}(t)) + \hat{B} u(t) \\[1ex]
\hat{w}(t) & = & \hat{C} \hat{v}(t) \\
\end{array}
\end{equation}
with 
$\hat{v} = (\hat{v}_{1}^\top,\ldots,\hat{v}_{m}^\top)^\top \in \real^{mn}$
and 
$\hat{w} = (\hat{w}_{1},\ldots,\hat{w}_{m})^\top \in \real^m$.
To define the involved matrices, we employ auxiliary arrays: 
the symmetric matrix $S(p) := (\Phi_i(p) \Phi_j(p)) \in \real^{m \times m}$ 
and the column vector $s(p) := (\Phi_i(p)) \in \real^m$.
Now the matrices of the linear parts in~(\ref{galerkin}) read as
$$ \begin{array}{ll} 
\hat{A} = \mathbb{E} [ S(\cdot) \otimes A(\cdot) ] 
\in \real^{mn \times mn}, \quad &
\hat{B} = \mathbb{E} [ s(\cdot) \otimes B(\cdot) ] 
\in \real^{mn}, \\[0.5ex]
\hat{C} = \mathbb{E} [ S(\cdot) \otimes C(\cdot) ] 
\in \real^{m \times mn}, \quad &
\hat{E} = \mathbb{E} [ S(\cdot) \otimes E(\cdot) ] 
\in \real^{mn \times mn} , \\
\end{array} $$
using the expected value~(\ref{expectedvalue}) componentwise 
and Kronecker products.
The nonlinear part is given by 
$\hat{F} := (\hat{F}_1^\top,\ldots,\hat{F}_m^\top)^\top$ 
with $\hat{F}_i \in \real^n$ and
\begin{equation} \label{fgal}
\hat{F}_{i} (\hat{v}) :=
\mathbb{E} \left[ 
F \left( \sum_{j=1}^m \hat{v}_j \Phi_j(\cdot) , \cdot \right) \Phi_i(\cdot) 
\right]
\end{equation}
for $i = 1,\ldots,m$. 

The expected values~(\ref{fgal}) represent probabilistic integrals, 
which often cannot be calculated analytically. 
An exception are functions~$F$ in~(\ref{dae}) consisting of 
polynomials with low degrees. 
We require approximations by quadrature formulas for the case of 
general nonlinear functions~$F$ in~(\ref{dae}). 
A quadrature rule or a sampling method is determined by its 
nodes $\{ p^{(1)},\ldots,p^{(k)} \} \subset \Pi$ and its
weights $\{ \gamma_1 , \ldots , \gamma_k \} \subset \real$. 
The associated approximation becomes
\begin{equation} \label{fcn-quadr}
\hat{F}_{i} (\hat{v}) \approx
\sum_{\ell = 1}^k \gamma_{\ell} \,
F \left( \sum_{j=1}^m \hat{v}_j \Phi_j(p^{(\ell)}) , p^{(\ell)} \right) 
\Phi_i(p^{(\ell)}) 
\end{equation}
for $i=1,\ldots,m$. 
Thus the computational effort for one evaluation of~$\hat{F}$ 
comprises $k$ evaluations of~$F$, whereas the summation is negligible.

Concerning initial value problems of the system~(\ref{galerkin}), 
the initial condition~(\ref{ivp}) of the original system~(\ref{dae}) 
has to be expanded as well. 
Consistent initial values are required in the case of DAEs. 

\subsection{Stochastic collocation method}
Alternatively, we consider a stochastic collocation 
(non-intrusive method) employing a quad\-ra\-ture rule or 
a sampling scheme, see~\cite{pulch14,xiu-hesthaven}.
This approach is determined uniquely by the nodes $\{ p^{(1)},\ldots,p^{(k)} \}$
and the weights $\{ \gamma_1 , \ldots , \gamma_k \}$. 
The approximations of the coefficient functions
belonging to the QoI~(\ref{qoi-method}) become
\begin{equation} \label{output-quadrature} 
\hat{w}_i(t) = \sum_{\ell=1}^k \gamma_{\ell} \Phi_i(p^{(\ell)}) \, 
C(p^{(\ell)}) x(t,p^{(\ell)}) 
\end{equation}
for $i=1,\ldots,m$.
Therein, $C(p^{(\ell)})$ denotes the evaluations of the output matrix
from~(\ref{dae}) for $\ell = 1,\ldots,k$.
Thus $k$ separate initial value problems~(\ref{dae}),(\ref{ivp}) have 
to be solved for the different nodes,
 which yields the variables~$x$.  

We write this approach as a single auxiliary system, 
which exhibits the form~(\ref{galerkin}) again. 
Due to a different meaning of the inner variables, we reformulate
\begin{equation} \label{collocation}
\begin{array}{rcl}
\hat{E} \dot{\hat{x}}(t) & = & 
\hat{A} \hat{x}(t) + \hat{F}(\hat{x}(t)) + \hat{B} u(t) \\[1ex]
\hat{w}(t) & = & \hat{C} \hat{x}(t) \\
\end{array}
\end{equation}
with $\hat{x}(t) = (x(t,p^{(1)})^\top, \ldots , x(t,p^{(k)})^\top)^\top$.
Consequently, the matrices read as 
$$ \begin{array}{ll}
\hat{A} = \begin{pmatrix}
A(p^{(1)}) & & \\
& \ddots & \\
& & A(p^{(k)}) \\
\end{pmatrix} \in \real^{kn \times kn} , \qquad & 
\hat{B} = \begin{pmatrix} 
B(p^{(1)}) \\ \vdots \\ B(p^{(k)}) \\
\end{pmatrix} \in \real^{kn} , \\[5ex]
$$ \hat{E} = \begin{pmatrix}
E(p^{(1)}) & & \\
& \ddots & \\
& & E(p^{(k)}) \\
\end{pmatrix} \in \real^{kn \times kn} . & \\
\end{array} $$
The matrix~$\hat{C} \in \real^{m \times kn}$ 
follows from the formula~(\ref{output-quadrature}).
The meaning of the outputs is the same as in the 
stochastic Galerkin system~(\ref{galerkin}). 
The nonlinear part becomes
$$ \hat{F}(\hat{x}(t)) = 
\begin{pmatrix} 
F(\hat{x}(t,p^{(1)}),p^{(1)}) \\ \vdots \\ F(\hat{x}(t,p^{(k)}),p^{(k)}) \\
\end{pmatrix} \in \real^{kn} . $$
The initial values follow from~(\ref{ivp}), 
i.e., $x(t_0,p^{(\ell)}) = x_0(p^{(\ell)})$ for $\ell = 1,\ldots,k$.
This auxiliary system was already derived and applied for 
linear dynamical systems in~\cite{pulch-scee2014,pulch17}. 
A similar single system was constructed 
in the case of a non-intrusive approach 
for dynamical systems with both stochastic noise and random parameters 
in~\cite{navarro}.

In the stochastic collocation method,
the system~(\ref{collocation}) is large and weakly coupled,
because it consists of the separate 
original systems~(\ref{dae}), which are connected just by the 
supply of the same input and the definition of the outputs. 



\section{Model order reduction}
\label{sec:mor}
Although the original dynamical system~(\ref{dae}) may be small or 
medium-sized, the dynamical systems~(\ref{galerkin}) and~(\ref{collocation}) 
from the numerical methods are large. 
Thus they represent excellent candidates for an MOR. 
Linear stochastic Galerkin systems were reduced successfully 
in~\cite{mi,pulch-scee2014,pulch-maten,pulch-maten-augustin,zou}.

\subsection{Projection-based model order reduction}
\label{sec:projection-mor}
The nonlinear dynamical systems~(\ref{galerkin}) and~(\ref{collocation}) 
represent non-parametric formulations. 
Without loss of generality, we consider the system~(\ref{galerkin}) 
as full-order model (FOM).
MOR yields a reduced-order model (ROM) in form of a 
smaller nonlinear dynamical system 
\begin{equation} \label{galerkin-reduced}
\begin{array}{rcl}
\bar{E} \dot{\bar{v}}(t) & = & 
\bar{A} \bar{v}(t) + \bar{F}(\bar{v}(t)) + \bar{B} u(t) \\[1ex]
\bar{w}(t) & = & \bar{C} \bar{v}(t) \\
\end{array}
\end{equation}
with matrices $\bar{A},\bar{E} \in \real^{r \times r}$, 
$\bar{B} \in \real^{r}$, $\bar{C} \in \real^{m \times r}$ 
and a function $\bar{F} : \real^r \rightarrow \real^r$. 
The state variables or inner variables are 
$\bar{v} : [t_0,t_{\rm end}] \rightarrow \real^r$, 
which require initial conditions $\bar{v}(t_0) = \bar{v}_0$. 
The input of~(\ref{galerkin-reduced}) coincides with the input 
of~(\ref{galerkin}). 
The aim is to achieve outputs satisfying $\bar{w}(t) \approx \hat{w}(t)$ 
for all~$t$.

In projection-based MOR, see~\cite{antoulas,freund}, 
projection matrices $T_{\rm l},T_{\rm r} \in \real^{mn \times r}$ 
of full rank are used. 
It follows that
\begin{equation} \label{mor-projections} 
\bar{A} = T_{\rm l}^\top \hat{A} T_{\rm r} , \quad 
\bar{B} = T_{\rm l}^\top \hat{B} , \quad 
\bar{C} = \hat{C} T_{\rm r} , \quad
\bar{E} = T_{\rm l}^\top \hat{A} T_{\rm r} , \quad 
\bar{F}(\bar{v}) = T_{\rm l}^\top \hat{F} ( T_{\rm r} \bar{v} ) . 
\end{equation}
A transient simulation of the dynamical systems often requires 
implicit integration schemes.
The smaller dimensionality of the system~(\ref{galerkin-reduced}) 
causes a lower computation work in the linear algebra algorithms.
However, the definition of the reduced nonlinear 
function~(\ref{mor-projections}) implies that 
a function evaluation of~$\bar{F}$ still takes an evaluation 
of the original function~$\hat{F}$. 
Significant reductions of the computational effort can be achieved only 
if one replaces the function evaluations of~$\hat{F}$ 
by cheaper approximations, 
which is also called hyper-reduction in this context. 
Several techniques have been proposed for this purpose 
like (discrete) empirical interpolation, missing point estimation 
or piecewise linear approaches. 
We refer to~\cite{chaturantabut} and the references therein. 
A hyper-reduction is not within the scope of this paper. 
Our main aim is the identification of sufficiently accurate 
low-dimensional representations.  

Furthermore, systems of DAEs often require more sophisticated MOR methods 
in comparison to systems of ODEs, see~\cite{benner-stykel,gugercin-stykel}.
Indeed the situation becomes more critical for higher-index DAEs 
(index larger than one).

\subsection{Proper orthogonal decomposition}
\label{sec:pod}
The projection-based ROM~(\ref{galerkin-reduced}),(\ref{mor-projections}) 
is determined uniquely by the matrices~$T_{\rm l}$ and $T_{\rm r}$.
We apply the POD to identify a projection matrix. 
A numerical integration scheme yields a solution of  
an initial value problem for the larger system~(\ref{galerkin}) 
or~(\ref{collocation}) in $[t_0,t_{\rm end}]$. 
In grid points $t_0 < t_1 < \cdots < t_{\ell-1} = t_{\rm end}$, 
we obtain approximations, which are called snapshots in this context. 
We collect them in a matrix 
\begin{equation} \label{pod_v}
V = \left( \hat{v}(t_0),\hat{v}(t_1),\ldots,\hat{v}(t_{\ell -1}) \right) 
\in \real^{mn \times \ell} . 
\end{equation}
Often it holds that $\ell \ll mn$.
Now a singular value decomposition (SVD), 
see~\cite[p.~76]{golub-loan} 
yields
\begin{equation} \label{svd}
V = U \begin{pmatrix} \Sigma \\ 0 \\ \end{pmatrix} Q^\top 
\end{equation}
with orthogonal matrices $U \in \real^{mn \times mn}$, 
$Q \in \real^{\ell \times \ell}$ and a diagonal matrix 
$\Sigma \in \real^{\ell \times \ell}$ 
with the singular values
$\sigma_1 \ge \sigma_2 \ge \cdots \ge \sigma_{\ell} \ge 0$.
Let $u_1, \ldots , u_{mn}$ be the columns of~$U$. 
The right-hand projection matrix is defined as
$$ T_{\rm r} = (u_1,u_2,\ldots,u_r) \in \real^{mn \times r} $$
for any $r \le \ell$. 
Thus just the singular vectors associated with the $r$ dominant 
singular values are included in the MOR scheme. 
The left-hand projection matrix is chosen by a Galerkin-type approach 
as $T_{\rm l} = T_{\rm r}$.

The POD technique requires a solution of an initial value problem 
of the FOM to generate the snapshots. 
Thus saving computational effort is achieved only if the 
ROM can be reused for other numerical simulations. 
These additional numerical simulations may involve
\begin{itemize}
\item other input signals~$u(t)$, 
\item other initial values~$\bar{v}(t_0)$, 
\item longer time intervals $[t_0,t']$ with $t' > t_{\rm end}$.
\end{itemize}
The third case will be examined for a test example 
in Section~\ref{sec:example}.

\subsection{Low-dimensional subspaces and approximation} 
The following derivation is applicable in the case of any 
projection-based MOR as discussed in Section~\ref{sec:projection-mor}.
The output of the ROM~(\ref{galerkin-reduced}) produces an approximation 
of the QoI in the original system~(\ref{dae}) by 
\begin{equation} \label{qoi-appr}
\bar{y}(t,p) = \displaystyle 
\sum_{i=1}^m \bar{w}_i(t) \Phi_i(p) 
\end{equation} 
for $t \in [t_0,t_{\rm end}]$ and $p \in \Pi$.
The $\ltwo$-error of this reduction with respect to~(\ref{qoi-method})
reads as
\begin{equation} \label{ltwo-error}
  \left\| \bar{y}(t,\cdot) - \hat{y}(t,\cdot) \right\|_{\ltwo} =
  \sqrt{ \sum_{i=1}^m \left( \bar{w}_i(t) - \hat{w}_i(t) \right)^2 }
  \qquad \mbox{for each}\;\; t ,
\end{equation}
due to the orthonormality of the basis and Parseval's equality.
Since a linear dependence of the QoI~$y$ on the state variables or 
inner variables~$x$ is assumed in~(\ref{dae}), 
the same derivation applies as in~\cite{pulch-matcom18}.
In~(\ref{galerkin-reduced}), the part $\bar{w} = \bar{C} \bar{v}$
with the output matrix $\bar{C} = (\bar{c}_{ij})$ yields
$$ \bar{w}_i(t) = \sum_{j=1}^{r} \bar{c}_{ij} \bar{v}_j(t)
\qquad \mbox{for} \;\; i=1,\ldots,m . $$

The approximation~(\ref{qoi-appr}) exhibits the formulation
$$ \bar{y}(t,p) =
\sum_{i=1}^m \left[ \sum_{j=1}^{r} \bar{c}_{ij} \bar{v}_j(t) \right] 
\Phi_i(p) =
\displaystyle \sum_{j=1}^{r} \bar{v}_{j}(t)
\left[ \sum_{i=1}^m \bar{c}_{ij} \Phi_i(p) \right] . $$
Thus new basis functions 
\begin{equation} \label{newbasis}
\Psi_j(p) := \sum_{i=1}^m \bar{c}_{ij} \Phi_i(p) 
\qquad \mbox{for}\;\; j=1,\ldots,r 
\end{equation}
are defined. 
We obtain the low-dimensional representation~(\ref{qoi-new}) 
with $\bar{w}_j = \bar{v}_j$ for $j=1,\ldots,r$.
Hence the coefficients are already identified by the MOR scheme. 

Let $\Phi = (\Phi_1,\ldots,\Phi_m)^\top$
and $\Psi = (\Psi_1,\ldots,\Psi_r)^\top$.
It holds that $\Psi = \bar{C}^\top \Phi$.
The set of functions~(\ref{newbasis}) is linearly independent 
in most of the cases. 
However, they are not orthogonal in general.
An orthonormal basis can be constructed from the original basis
by an SVD of the output matrix
\begin{equation} \label{svd-output}
\bar{C} = U \begin{pmatrix} \Sigma \\ 0 \\ \end{pmatrix} Q^\top 
\end{equation}
including orthogonal matrices $U \in \real^{m \times m}$, 
$Q \in \real^{r \times r}$ and a diagonal matrix 
$\Sigma \in \real^{r \times r}$ 
with the singular values
$\sigma_1 \ge \sigma_2 \ge \cdots \ge \sigma_{r} \ge 0$.
Thus we obtain the alternative basis functions
\begin{equation} \label{newbasis2}
  \Psi_j^*(p) = \sum_{i=1}^m u_{ij} \Phi_i(p)
\end{equation}
with coefficients from the matrix~$U$.
If follows that $\{ \Psi_1^* , \ldots , \Psi_r^* \}$ is an
orthonormal basis spanning the same subspace as
$\{ \Psi_1 , \ldots , \Psi_r \}$.
The alternative basis can be deflated to
$\{ \Psi_1^* , \ldots , \Psi_{k}^* \}$ for a $k \le r$
to remove a (numerical) rank deficiency. 
The details are explained in~\cite{pulch-matcom18}.

The original orthonormal basis $\{ \Phi_1,\ldots,\Phi_m \}$ owns the
advantage that the first coefficient yields the expected value
due to $\Phi_1 \equiv 1$.
The first and second moment can be obtained as well in the novel
orthonormal basis.
Expected value and variance read as
$$ \begin{array}{rcl}
  \displaystyle
  \mathbb{E} \left[ \sum_{j=1}^r w_j(t) \Psi_j^*(\cdot) \right]
  & = & \displaystyle \sum_{j=1}^r w_j(t) u_{1j} \\
  \displaystyle
      {\rm Var}  \left[ \sum_{j=1}^r w_j(t) \Psi_j^*(\cdot) \right]
  & = & \displaystyle \left( \sum_{j=1}^r w_j(t)^2 \right) -
  \left( \sum_{j=1}^r w_j(t) u_{1j} \right)^2 \\
\end{array} $$
due to~(\ref{newbasis2}).
Hence just the SVD~(\ref{svd-output}) has to be computed at the beginning.
The moments do not require significant additional work.
Yet the only advantage of the basis 
$\{ \Psi_1^* , \ldots , \Psi_r^* \}$ in comparison to
$\{ \Phi_1,\ldots,\Phi_m \}$ is the reduced dimensionality.

\subsection{Best approximation}
A projection-based MOR~(\ref{galerkin-reduced}),(\ref{mor-projections}) 
identifies a low-dimensional representation (\ref{qoi-new}) 
including its time-dependent coefficient functions. 
In addition, we obtain a best approximation within the subspace spanned by
the new basis~(\ref{newbasis}), where just information from the 
right-hand projection matrix $T_{\rm r}$ is used.

The best approximation 
\begin{equation} \label{bestappr}
y^*(t,p) = \sum_{j=1}^r w_j^*(t) \Psi_j(p) 
\end{equation}
with $w^* = (w_1^*,\ldots,w_r^*)^\top \in \real^r$
is defined by the optimisation problem
\begin{equation} \label{bestappr2} 
\min_{w(t)} \;
\left\| \hat{y}(t,\cdot) - \sum_{j=1}^r w_j(t) \Psi_j(\cdot) 
\right\|_{\ltwo} 
\end{equation}
with $w = (w_1,\ldots,w_r)^\top \in \real^r$ 
pointwise for $t \ge t_0$. 
This best approximation can be computed from the solution of the FOM 
by a linear least squares problem.

\clearpage

\begin{theorem}
Let $\hat{w} = (\hat{w}_1,\ldots,\hat{w}_m)^\top$ be the output 
of the dynamical system~(\ref{galerkin}) or~(\ref{collocation}) 
for an initial value problem. 
Let the output matrix $\bar{C} = \hat{C} T_{\rm r}$ have full rank.
The coefficients $w^* = (w_1^*,\ldots,w_r^*)$ of the 
best approximation~(\ref{bestappr}) minimising~(\ref{bestappr2}) 
are the solution of the linear least squares problem
\begin{equation} \label{leastsquares}
\min_{w(t)} 
\left\| \hat{w}(t) - \bar{C} w(t) \right\|_2 
\end{equation}
pointwise for~$t \ge t_0$.
\end{theorem}

Proof:

We omit the dependence on time for notational convenience. 
The Hilbert space norm can be expressed by inner products 
$$ \left\| \hat{y} - \bar{y} \right\|_{\ltwo}^2 = \;
< \hat{y} , \hat{y} > - 2 < \hat{y} , \bar{y} > + < \bar{y} , \bar{y} > . $$
We obtain
$$  < \hat{y} , \bar{y} > = 
\sum_{i=1}^m \sum_{j=1}^r \hat{w}_i \bar{w}_j < \Phi_i , \Psi_j > 
\quad \mbox{and} \quad
< \bar{y} , \bar{y} > = 
\sum_{i,j=1}^r \bar{w}_i \bar{w}_j < \Psi_i , \Psi_j > . $$
Using~(\ref{newbasis}), basic calculations yield
$$ < \Phi_i , \Psi_j > = \bar{c}_{ij} 
\quad \mbox{and} \quad
< \Psi_i , \Psi_j > = \sum_{k=1}^m \bar{c}_{ki} \bar{c}_{kj} $$
with $\bar{C} = (\bar{c}_{ij})$. 
It follows that
$$ \left\| \hat{y} - \bar{y} \right\|_{\ltwo}^2 = \;
< \hat{y} , \hat{y} > - 2 \bar{w}^\top \bar{C}^\top \hat{w} 
+ \bar{w}^\top \bar{C}^\top \bar{C} \bar{w} . $$
The degrees of freedom are $\bar{w} \in \real^r$, 
whereas $\hat{w} \in \real^m$ is constant. 
A necessary condition for a minimum is a vanishing gradient. 
Thus we achieve
$$ 2 \bar{C}^\top \bar{C} \bar{w} - 2 \bar{C}^\top \hat{w} = 0 , $$
which is equivalent to the normal equation, 
see~\cite[p.~232]{stoer-bulirsch}, 
associated with the linear least squares problem~(\ref{leastsquares}).
\hfill $\Box$
\bigskip

The least squares problem~(\ref{leastsquares}) implies the formula
\begin{equation} \label{bestappr3}
w^*(t) = ( \bar{C}^\top \bar{C} )^{-1} \bar{C}^\top \hat{w} (t) 
\qquad \mbox{for all}\;\; t 
\end{equation}
provided that the output matrix exhibits full rank.
Due to~(\ref{bestappr3}), the best approximation~(\ref{bestappr}) 
is continuous in time, because the solution of the 
dynamical system~(\ref{galerkin}) or~(\ref{collocation}) is 
assumed to be continuous. 
The computation of the best approximation requires to solve the FOM. 
However, the application of the formula~(\ref{bestappr3}) afterwards 
is cheap, because the transformation matrix is time-invariant. 
A QR-decomposition of the output matrix can be reused at all time points.

\subsection{Error analysis}
\label{sec:erroranalysis}
The difference $\hat{y} - \bar{y}$ between the QoI~(\ref{qoi-method}) 
from a numerical method and the approximation~(\ref{qoi-new}) 
is exactly the error of the MOR approach. 
Hence this error depends on the individual choice of the 
projection-based MOR method.

A more detailed examination is feasible for the best approximation 
using the POD method of Section~\ref{sec:pod}.  
Again the Galerkin system~(\ref{galerkin}) is considered 
without loss of generality, since 
the analysis also applies to the collocation system~(\ref{collocation}).
We obtain the following property of the approximation quality with respect to 
the state variables or the inner variables.

\begin{lemma} \label{lemma1}
Let $\hat{v} \in C^1$ be the exact solution of an initial value problem 
of the dynamical system~(\ref{galerkin}).  
Let $T_{\rm r} \in \real^{mn \times r}$ be the projection matrix from  
the POD approach with $\ell > r$ snapshots.
For each $t \in [t_0,t_{\rm end}]$, the solution~$v^*(t) \in \real^r$ 
of the linear least squares problem
\begin{equation} \label{leastsquares-v}
\min_{v(t)} \; \left\| \hat{v}(t) - T_{\rm r} v(t) \right\|_2 
\end{equation}
satisfies the estimate
$$ \left\| \hat{v}(t) - T_{\rm r} v^*(t) \right\|_2 \le 
\sigma_{r+1} + 
\sqrt{mn} \; \Delta t \max_{t_0 \le \tau \le t_{\rm end}} 
\left\| \dot{\hat{v}} (\tau) \right\|_{\infty} $$
with the time step size 
$\Delta t = \max \{ t_{j} - t_{j-1} \; : \; j=1,\ldots,\ell-1 \}$ 
and the singular value $\sigma_{r+1}$ from the POD.
\end{lemma}

Proof:

The intermediate value theorem of differential calculus yields componentwise 
$\hat{v}_i(t) = \hat{v}_i(t_j) + (t-t_j) \dot{\hat{v}}_i(\xi)$ 
for each $t \in [t_j,t_{j+1}]$ and $i=1,\ldots,mn$ 
with intermediate values~$\xi$.
It follows that
\begin{equation} \label{bound-taylor}
\left\| \hat{v}(t) - \hat{v}(t_j) \right\|_{\infty} \le 
\Delta t \max_{t_0 \le \tau \le t_{\rm end}} 
\left\| \dot{\hat{v}} (\tau) \right\|_{\infty} 
\quad \mbox{for} \;\; t \in [t_j,t_{j+1}] .
\end{equation}
Since $v^*$ represents an optimum, we obtain
$$ \left\| \hat{v}(t) - T_{\rm r} v^*(t) \right\|_2 \le 
\left\| \hat{v}(t) - T_{\rm r} \breve{v} \right\|_2 \le 
\left\| \hat{v}(t) -  \hat{v}(t_j) \right\|_2 + 
\left\| \hat{v}(t_j) - T_{\rm r} \breve{v} \right\|_2 $$
for any $\breve{v} \in \real^r$. 
The first term can be bounded by 
$$ \left\| \hat{v}(t) -  \hat{v}(t_j) \right\|_2 \le 
\sqrt{mn} \left\| \hat{v}(t) -  \hat{v}(t_j) \right\|_{\infty} $$ 
and the estimate~(\ref{bound-taylor}) for each $t \in [t_j,t_{j+1}]$.
For the second term, we apply the equality $\hat{v}(t_j) = V e_j$ with 
the matrix~(\ref{pod_v}) and a canonical unit vector $e_j \in \real^{\ell}$.
Let $V_r \in \real^{mn \times \ell}$ be the closest rank-$r$ approximation  
identified by the SVD, see~\cite[p.~79]{golub-loan}.
It holds that 
$$ V_r = U \begin{pmatrix} \Sigma_r \\ 0 \\ \end{pmatrix} Q^\top
= \sum_{i=1}^r \sigma_i u_i q_i^\top , $$
where the diagonal matrix $\Sigma_r \in \real^{\ell \times \ell}$ 
contains only the $r$ dominant singular values. 
Furthermore, we obtain $V_r = T_{\rm r} \widetilde{\Sigma}_r Q^\top$ 
with a modified diagonal matrix~$\widetilde{\Sigma}_r$,
 because $T_{\rm r}$ consists of the first~$r$ columns of~$U$. 
Thus we choose the vector $\breve{v} = \widetilde{\Sigma}_r Q^\top e_j$. 
It follows that
$$ \left\| \hat{v}(t_j) - T_{\rm r} \breve{v} \right\|_2 = 
\left\| V e_j - V_r e_j \right\|_2 \le 
\| V - V_r \|_2 \| e_j \|_2 = 
\| V - V_r \|_2 = \sigma_{r+1} $$
due to the error estimate for the spectral norm 
in~\cite[p.~79]{golub-loan}. 
\hfill $\Box$

\bigskip

The solutions of the least squares problems~(\ref{leastsquares-v}) 
are smooth again, because $C^1$-solutions of the FOM are assumed. 
However, we do not require this smoothness of the optimum in the following.
Now we show an error bound on the QoI.

\begin{theorem}
Let $\hat{v} \in C^1$ be the solution and $\hat{y}$ be the QoI 
of an initial value problem of the dynamical system~(\ref{galerkin}).
The best approximation~$y^*$ with respect to the subspace, 
which is identified by the POD with $\ell > r$ snapshots, 
satisfies the error bound
\begin{equation} \label{errorbound}
\left\| \hat{y}(t,\cdot) - y^*(t,\cdot) \right\|_{\mathcal{L}_2(\Pi,\rho)}
\le \| \hat{C} \|_2 
\left( \sigma_{r+1} + 
\sqrt{mn} \; \Delta t \max_{t_0 \le \tau \le t_{\rm end}} 
\left\| \dot{\hat{v}} (\tau) \right\|_{\infty} \right) 
\end{equation}
for each $t \in [t_0,t_{\rm end}]$
with the time step size 
$\Delta t = \max \{ t_{j} - t_{j-1} \; : \; j=1,\ldots,\ell-1 \}$ 
and the singular value $\sigma_{r+1}$ from the POD.
\end{theorem}

Proof:

It holds that $\Psi = \bar{C}^\top \Phi$ with 
$\Psi := ( \Psi_1,\ldots,\Psi_r )^\top$ and 
$\Phi := ( \Phi_1 , \ldots , \Phi_m )^\top$ 
due to~(\ref{newbasis}).
The best approximation exhibits a representation
$$ y^*(t,p) = \sum_{j=1}^r w_j^*(t) \Psi_j(p) = 
w^*(t)^\top \Psi(p) = ( \bar{C} w^*(t) )^\top \Phi(p) . $$
with $w^* = (w_1^*,\ldots,w_r^*)^\top$. 
The POD yields the projection matrix $T_{\rm r}$. 
The orthonormality of the basis functions~$\Phi$ 
allows for the application of Parseval's equality.
For each $t \in [t_0,t_{\rm end}]$, it follows that
$$ \begin{array}{rclcl}
\left\| \hat{y}(t,\cdot) - y^*(t,\cdot) \right\|_{\mathcal{L}_2(\Pi,\rho)}
& = & \| \hat{w}(t) - \bar{C} w^*(t) \|_2 
& \le & \| \hat{w}(t) - \bar{C} \breve{w}(t) \|_2 \\[2ex] 
& = & \| \hat{C} \hat{v}(t) - \hat{C} T_{\rm r} \breve{w}(t) \|_2 
& \le &  \| \hat{C} \|_2 
\| \hat{v}(t) - T_{\rm r} \breve{w}(t) \|_2   \\
\end{array} $$
for any $\breve{w} \in \real^r$. 
Thus Lemma~\ref{lemma1} yields the bound~(\ref{errorbound}) 
by the choice $\breve{w} = v^*$ from 
the least squares problem~(\ref{leastsquares-v}).
\hfill $\Box$

\bigskip

The spectral norm of the output matrix is often uncritical in 
the estimate~(\ref{errorbound}). 
For example, if the single output of~(\ref{dae}) is just a state variable 
or inner variable, then it follows that $\| \hat{C} \|_2 = 1$.
The error bound~(\ref{errorbound}) becomes small in the case of a 
fast decay of the singular values and a sufficiently small time step size. 
On the one hand, the singular values are computed a priori, 
which yields the first term of the estimate. 
On the other hand, the second term is not computable, because the 
involved derivatives are unknown. 


\section{Illustrative examples}
\label{sec:example}
Now we apply an MOR approach to the stochastic Galerkin system 
as well as a stochastic collocation formulation for two test examples.
All numerical calculations are performed within the software package 
MATLAB. 

\subsection{Scrapie model}
\label{sec:scrapie}
We use a model of the scrapie disease 
from~\cite[p.~37]{deuflhard-bornemann}.
Reaction kinetics yields an autonomous system of ODEs 
\begin{equation} \label{scrapie} 
\begin{array}{rcl}
\dot{x}_1 & = & - p_1 x_1 + p_2 x_2 - p_5 x_1 x_3 \\
\dot{x}_2 & = & p_1 x_1 - p_2 x_2 - 2 p_3 x_2^2 
+ 2 p_4 x_3 + p_5 x_1 x_3 \\
\dot{x}_3 & = & p_3 x_2^2 - p_4 x_3 \\
\end{array} 
\end{equation}
for three molecular concentrations $x_1,x_2,x_3$. 
The system~(\ref{scrapie}) owns the form~(\ref{dae}) with 
$B=0$ and an identity matrix~$E$. 
As output, we examine the three concentrations separately. 
The five parameters are reaction constants, where 
we arrange the nominal values $p^* = (10^{-5},0.1,1,10^{-4},0.1)^\top$ 
from~\cite{deuflhard-bornemann}. 
Initial values~(\ref{ivp}) are chosen as 
$x(0,p) = (1,0,0.1)^\top$ for all~$p$.  
The time interval $[0,500]$ is considered.

In the stochastic modelling, we replace the parameters by independent 
uniformly distributed random variables with 10\% variation around 
the nominal values~$p^*$. 
The PC expansion~(\ref{pce}) includes multivariate polynomials, 
which are the products of univariate Legendre polynomials.
Our truncated expansion~(\ref{pce-appr}) involves all polynomials up to 
total degree three, which implies $m=56$ basis functions. 
The stochastic Galerkin method yields a dynamical system~(\ref{galerkin}) 
with $mn = 168$ state variables. 
Since just quadratic nonlinearities are given in the right-hand side 
of the system~(\ref{scrapie}), the right-hand sides 
of~(\ref{galerkin}) can be evaluated exactly except for roundoff errors. 
In a stochastic collocation approach, we consider a tensor-product 
formula of the Gauss-Legendre quadrature, 
see~\cite[p.~171]{stoer-bulirsch}, with $k=3^5=243$ nodes. 
The accompanying dynamical system~(\ref{collocation}) exhibits 
$kn = 729$ state variables. 

\begin{figure}[t]
\begin{center}
\includegraphics[width=6.5cm]{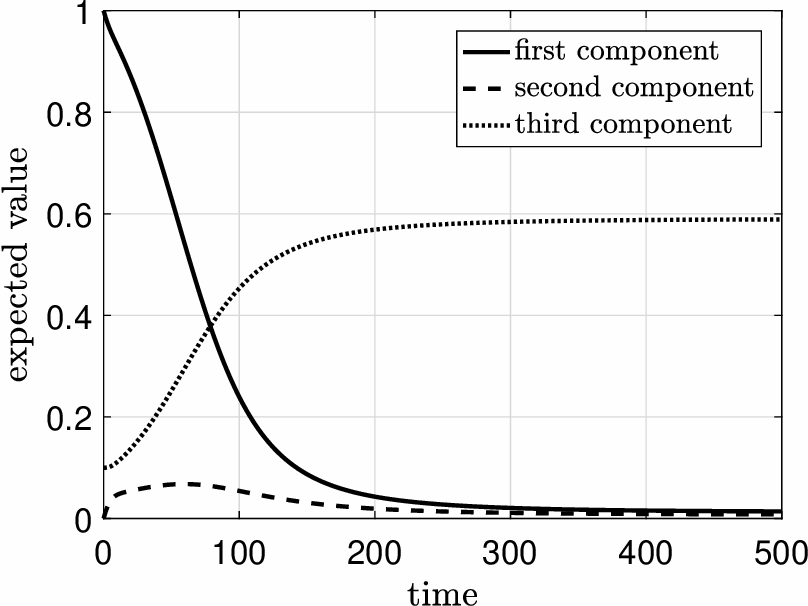}
\hspace{5mm}
\includegraphics[width=6.5cm]{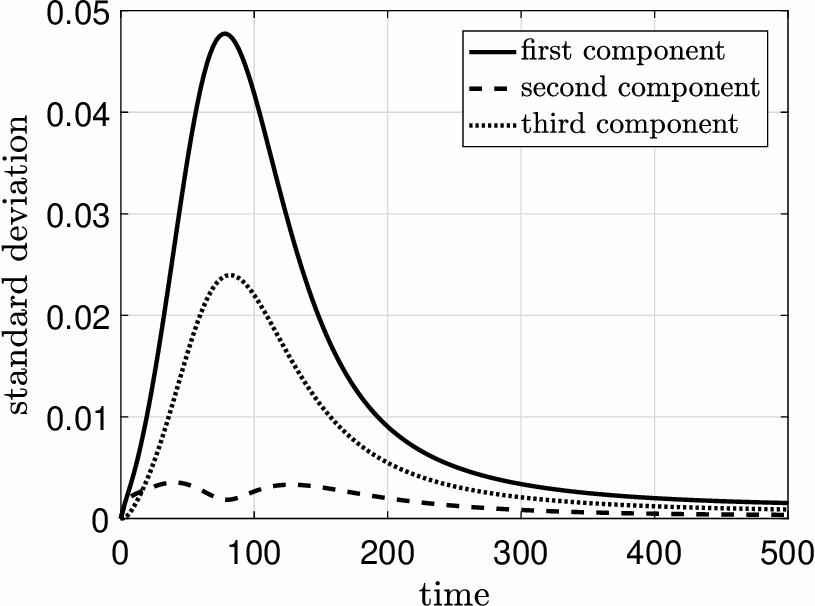}
\end{center}
\caption{Expected values (left) and standard deviations (right) for 
the random concentrations in the scrapie model.}
\label{fig:scrapie_stats}
\end{figure}

\begin{figure}[t]
\begin{center}
\includegraphics[width=6.5cm]{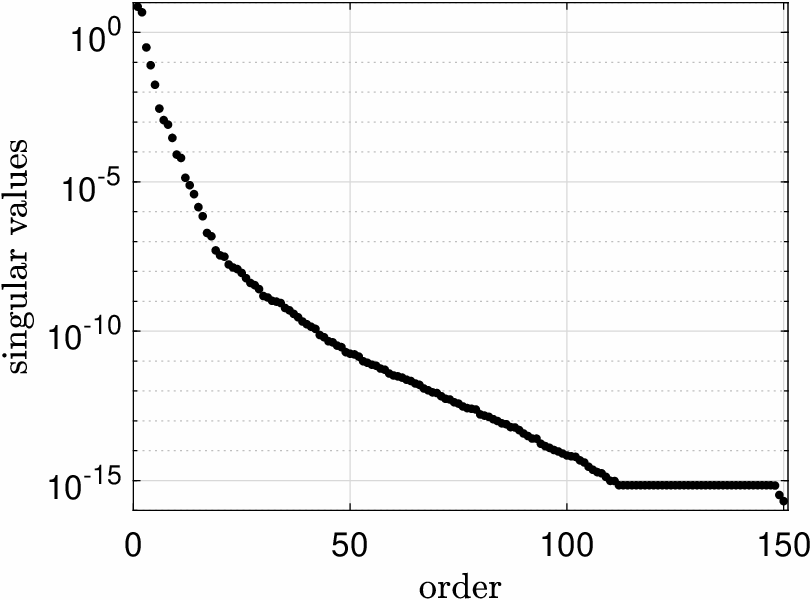}
\hspace{5mm}
\includegraphics[width=6.5cm]{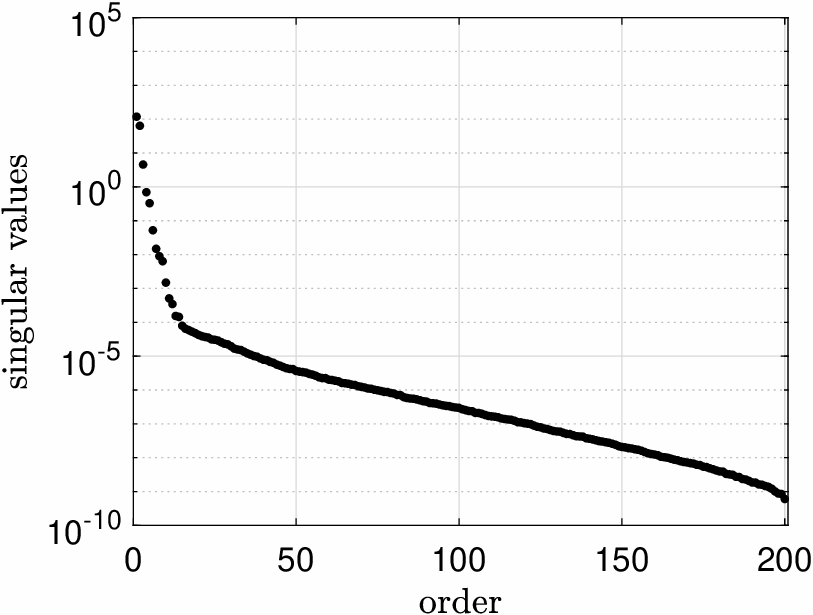}
\end{center}
\caption{Singular values from the POD in the stochastic Galerkin method 
(left) and the stochastic collocation method (right) 
for the scrapie example.}
\label{fig:scrapie_pod}
\end{figure}

The trapezoidal rule yields the numerical solutions of all 
initial value problems in this test example. 
Variable time step sizes are determined by a local error control 
with relative tolerance~$\varepsilon_{\rm rel}$ and absolute 
tolerance~$\varepsilon_{\rm abs}$. 
In the stochastic Galerkin method, we apply the choice 
$\varepsilon_{\rm rel} = 10^{-4}$, $\varepsilon_{\rm abs} = 10^{-6}$ 
for the computation of the snapshots, which causes 149 steps 
and thus 150 snapshots including the initial values. 
These snapshots are also used to obtain an approximation of the 
expected values  (first coefficient)  as well as
the standard deviations  (other coefficients) 
of the three random concentrations shown in Figure~\ref{fig:scrapie_stats}. 
Now the POD method requires the SVD~(\ref{svd}). 
Figure~\ref{fig:scrapie_pod} (left) illustrates the computed 
singular values, which decay rapidly. 
The projection matrices $T_{\rm l},T_{\rm r}$ and the 
ROM~(\ref{galerkin-reduced}) follow from the 
POD technique for user-defined reduced dimensions.
Since the dimensionality of the FOM~(\ref{galerkin}) is relatively small, 
we cannot expect saving computational effort by an MOR. 
Nevertheless, we discuss the MOR to show the feasibility of the approach 
for the identification of a low-dimensional representation.
 
\begin{figure}[t]
\begin{center}
\includegraphics[width=6.5cm]{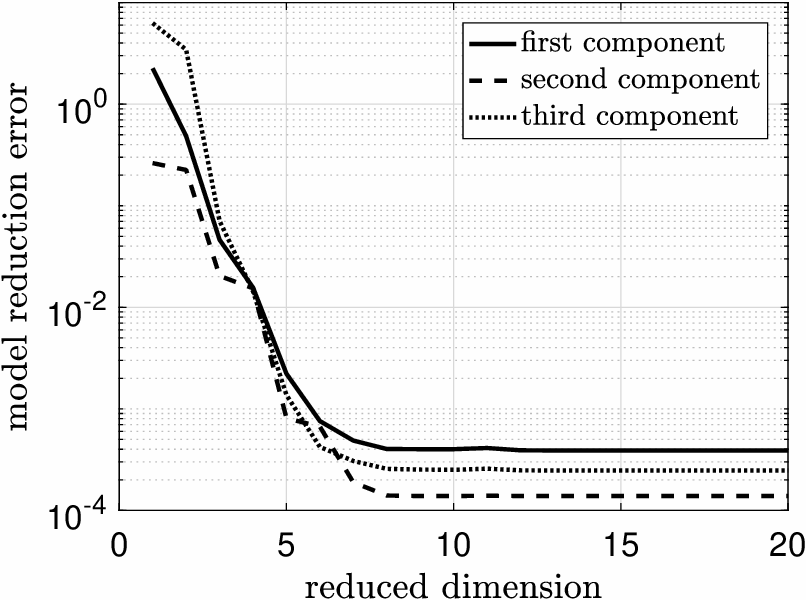}
\hspace{5mm}
\includegraphics[width=6.5cm]{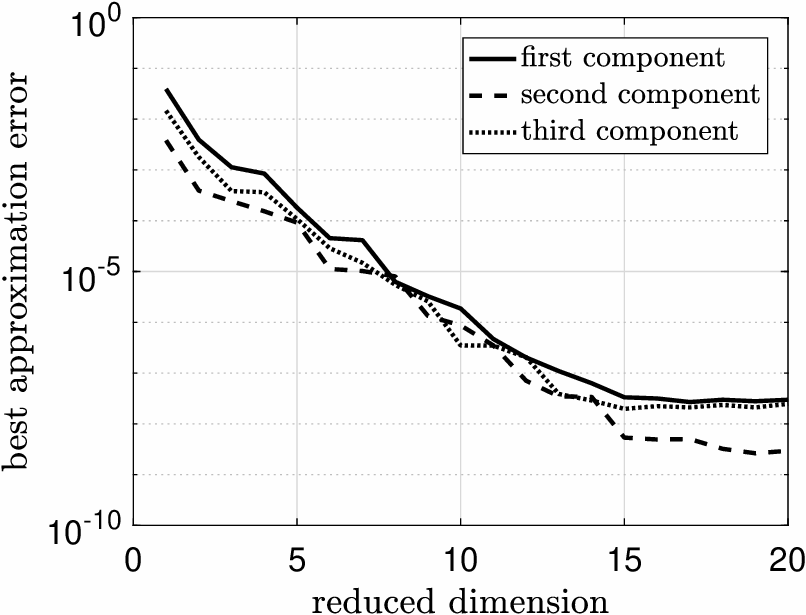}
\end{center}
\caption{Maximum $\ltwo$-errors for low-dimensional 
representations from MOR (left) and from best approximation (right) 
in Galerkin method for the scrapie model.}
\label{fig:scrapie_errors_gal}
\end{figure}

Now both the FOM~(\ref{galerkin}) and the 
ROMs~(\ref{galerkin-reduced}) are integrated
with accuracy requirements 
$\varepsilon_{\rm rel} = 10^{-3}$, $\varepsilon_{\rm abs} = 10^{-6}$. 
Initial values are transformed via $\bar{v}(0) = T_{\rm r}^\top \hat{v}(0)$.
We obtain the approximation~(\ref{qoi-appr}) from solving 
the ROM as well as the best approximation~(\ref{bestappr}) 
for each concentration as QoI. 
 The $\ltwo$-error~(\ref{ltwo-error}) between FOM and 
ROM can be determined pointwise in time. 
We evaluate the approximations at 200 equidistant time points 
in the total time interval using interpolation in time. 
 Figure~\ref{fig:scrapie_errors_gal} depicts the
$\ltwo$-error~(\ref{ltwo-error}) for different reduced dimensions,
where the maximum error of all time points is determined. 
We recognise that the MOR error decreases rapidly for increasing dimensions
at the beginning and then stagnates. 
The error becomes dominated by the quality of the snapshots, 
which depends on the density of the associated time points, 
and the error of the time integration. 
In contrast, the error of the best approximation decreases further for 
increasing dimensions and achieves a much lower magnitude.

\begin{figure}[t]
\begin{center}
\includegraphics[width=6.5cm]{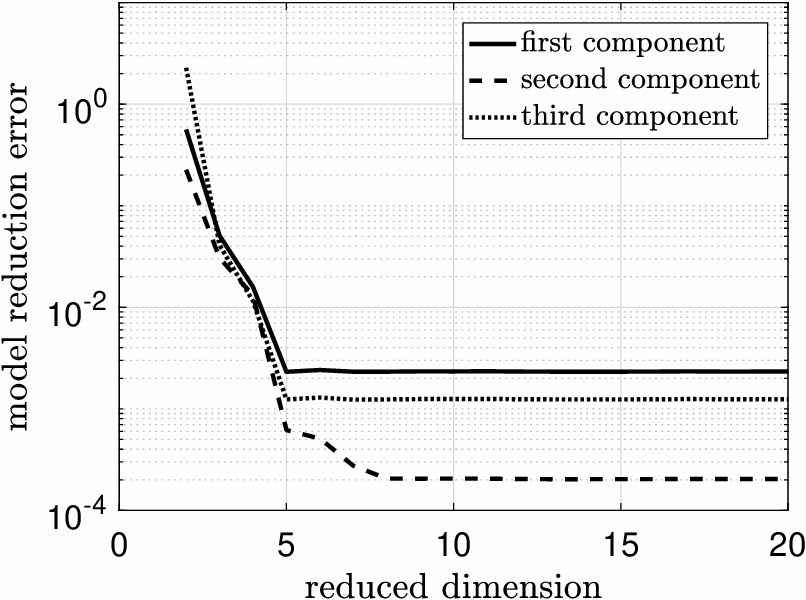}
\hspace{5mm}
\includegraphics[width=6.5cm]{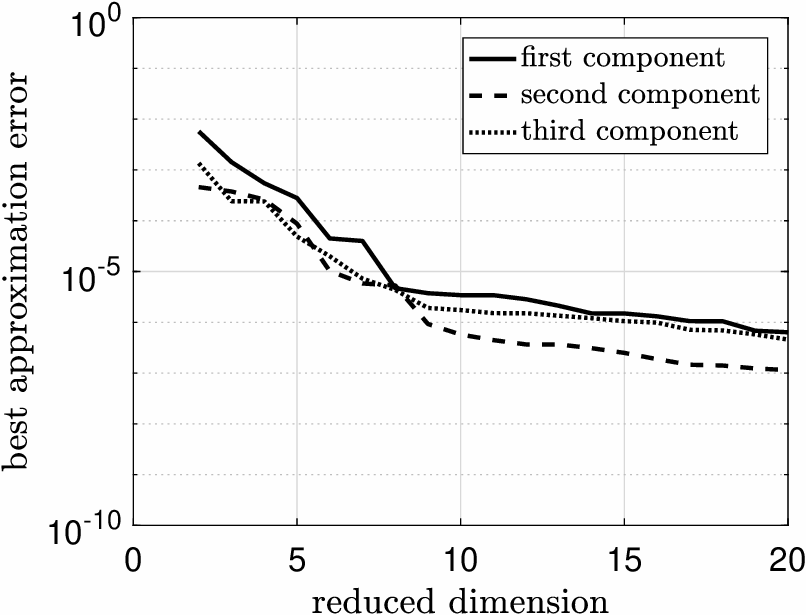}
\end{center}
\caption{Maximum $\ltwo$-errors for low-dimensional 
representations from MOR (left) and from best approximation (right) 
in collocation method for the scrapie model.}
\label{fig:scrapie_errors_col}
\end{figure}

We repeat the strategy for the stochastic collocation method now. 
Each system~(\ref{scrapie}) is integrated separately for the nodes 
of the quadrature with different step sizes for tolerances
$\varepsilon_{\rm rel} = 10^{-4}$, $\varepsilon_{\rm abs} = 10^{-6}$. 
Interpolation yields 200 snapshots of the weakly 
coupled system~(\ref{collocation}) at equidistant time points.
The result of the POD is shown in Figure~\ref{fig:scrapie_pod} (right), 
which demonstrates a fast decay of the singular values again.
We impose the accuracy requirements
$\varepsilon_{\rm rel} = 10^{-3}$, $\varepsilon_{\rm abs} = 10^{-6}$
on the integration of both the FOM~(\ref{collocation}) and the ROMs. 
Initial values are obtained by $\bar{x}(0) = T_{\rm r}^\top \hat{x}(0)$. 
Figure~\ref{fig:scrapie_errors_col} illustrates the maximum 
$\ltwo$-errors~(\ref{ltwo-error}) for the approximations~(\ref{qoi-appr}) 
and~(\ref{bestappr}) of the QoIs 
in the case of different reduced dimensions. 
A comparison to Figure~\ref{fig:scrapie_errors_gal} shows that 
the quality of the approximations is slightly worse in this 
stochastic collocation approach. 
We remark that the accuracy of the stochastic Galerkin method or
the stochastic collocation method with respect to the exact random QoI 
in~(\ref{dae}) is not compared, because this item is out of scope.
We determine the error of the MOR separately in the Galerkin technique
and the collocation scheme.

\begin{figure}[t]
\begin{center}
\includegraphics[width=8cm]{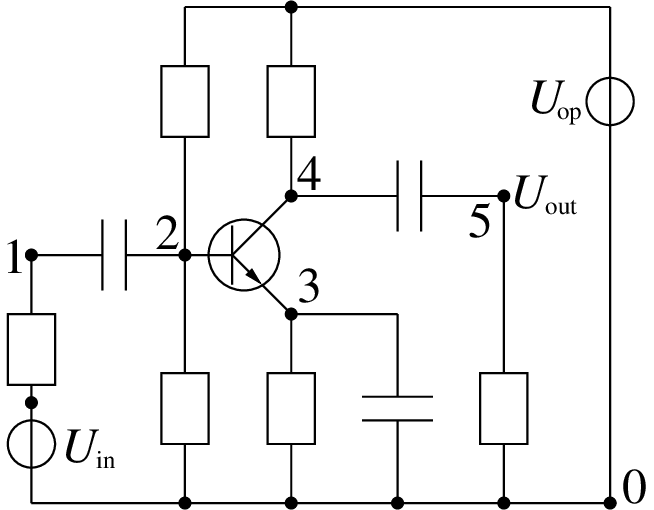} 
\end{center}
\caption{Electric circuit of a transistor amplifier.}
\label{fig:transistoramplifier}
\end{figure}

\subsection{Transistor amplifier circuit}
\label{sec:example-modelling}
\begin{figure}
\begin{center}
\includegraphics[width=8cm]{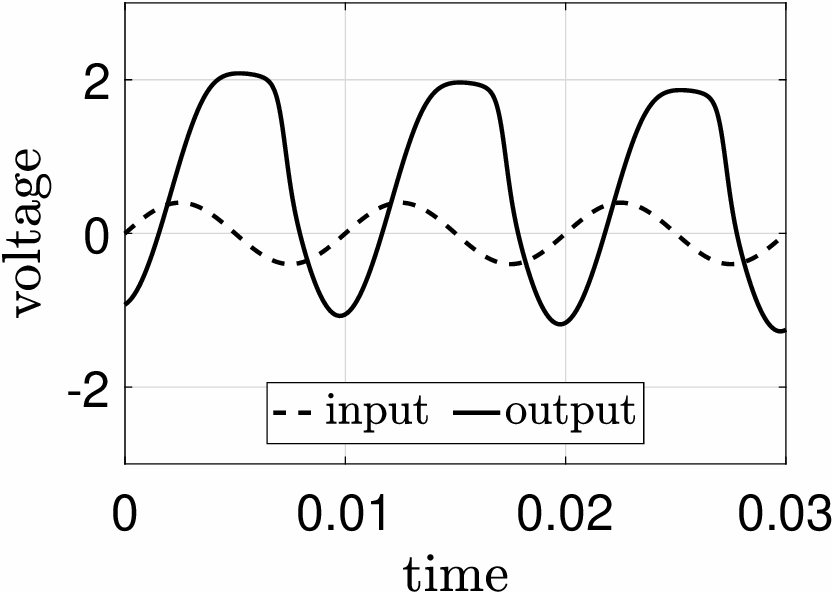}
\end{center}
\caption{Input voltage and output voltage of the transistor amplifier 
in the case of constant parameters.}
\label{fig:output}
\end{figure}

We consider the electric circuit of a transistor amplifier 
depicted in Figure~\ref{fig:transistoramplifier}. 
A mathematical modelling yields a dynamical system~(\ref{dae}) 
of DAEs ($n=5$) with differential index one, 
which is given in~\cite[p.~377]{hairer2}. 
The system is nonlinear, where the mapping~$F$ models 
the bipolar transistor of the circuit and thus involves exponential functions. 
Three capacitances and six resistances are included in the 
matrices~$E$ and~$A$, respectively. 
The matrix~$B \in \real^{5 \times 2}$ inserts the 
inputs $u = (U_{\rm in},U_{\rm op})^\top$. 
The operating voltage $U_{\rm op}$ is constant. 
We put the value~$U_{\rm op}$ into the matrix~$B$ and the 
input becomes $u = (U_{\rm in},\widetilde{U}_{\rm op})^\top$ 
with $\widetilde{U}_{\rm op} \equiv 1$, 
because $U_{\rm op}$ is a parameter in the matrix~$B$ now.
A single time-varying input is supplied by the voltage source $U_{\rm in}$, 
which we select as the sinusoidal input 
$$ U_{\rm in}(t) = 0.4 \sin \left( \textstyle \frac{2\pi}{T} t \right) 
\qquad \mbox{with} \;\; T = 0.01 . $$
The unknowns of the system consist of the five node voltages. 
The QoI is defined as the output voltage $U_{\rm out}$ at the fifth node.
Thus the matrix~$C$ in~(\ref{dae}) is just a unit vector.
Figure~\ref{fig:output} shows the numerical solution of an 
initial value problem~(\ref{dae}),(\ref{ivp}) for constant 
physical parameters from~\cite{hairer2}. 

\begin{figure}
\begin{center}
\includegraphics[width=6.5cm]{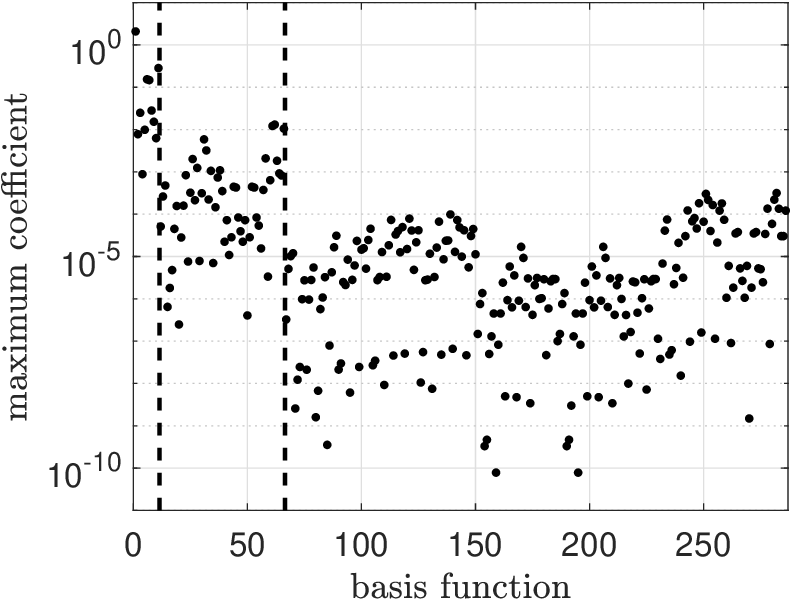}
\hspace{5mm}
\includegraphics[width=6.5cm]{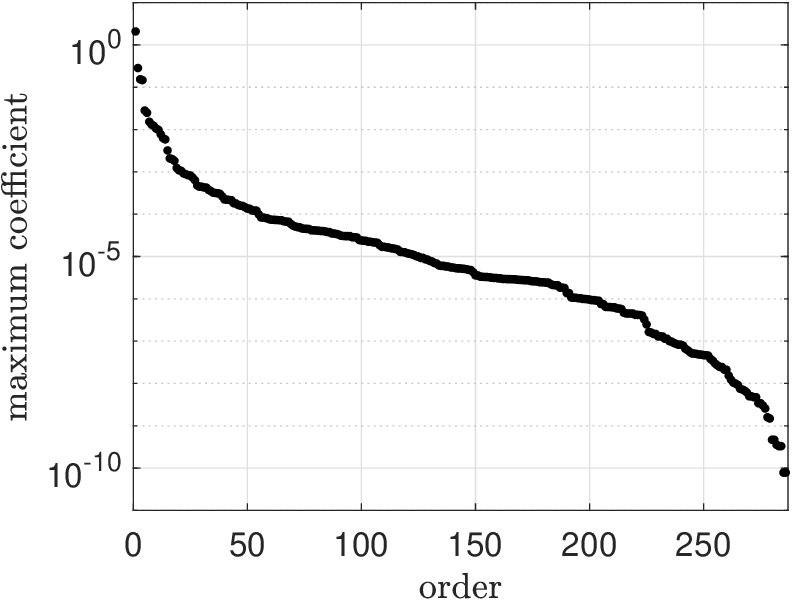}
\end{center}
\caption{Maximum coefficient functions (left) 
and their descending order (right) for the random output voltage 
in the transistor amplifier. 
(The dashes lines separate coefficients for polynomials of degree zero/one, 
two and three.)}
\label{fig:transamp_coeff}
\end{figure}

In this test example, all numerical solutions of initial value problems 
are computed by the backward differentiation formulas (BDF), 
see~\cite[p.~531]{stoer-bulirsch}. 
A local error control with tolerances
$\varepsilon_{\rm rel}$ and $\varepsilon_{\rm abs}$ 
yields adaptive time step sizes as well as adaptive orders (1~to~5). 

\begin{figure}
\begin{center}
\includegraphics[width=6.5cm]{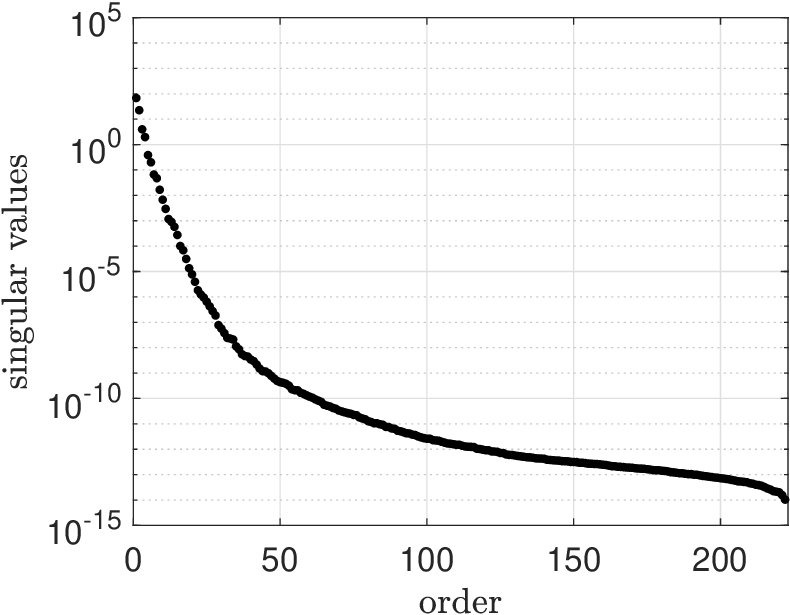}
\hspace{5mm}
\includegraphics[width=6.5cm]{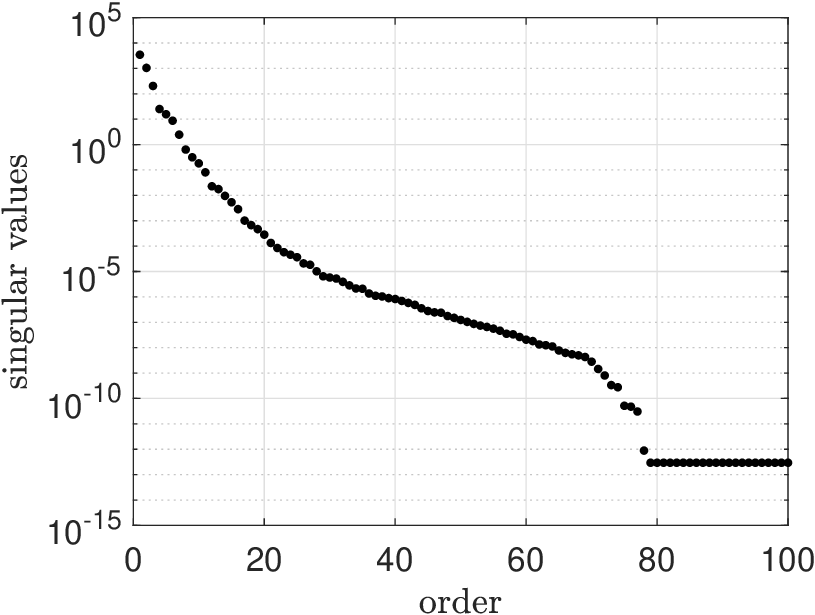}
\end{center}
\caption{Singular values from the POD in the stochastic Galerkin method 
(left) and the stochastic collocation method (right) 
for the transistor amplifier.}
\label{fig:transamp_pod}
\end{figure}

For an uncertainty quantification, we choose all capacitances, all resistances 
and the operating voltage as random parameters ($q=10$) with independent 
uniform distributions varying 1\% around the constant parameters from above.
Concerning the orthogonal expansion~(\ref{pce-appr}), 
we apply all basis polynomials up to total degree three 
and obtain $m=286$ terms.

The stochastic Galerkin method generates a DAE system~(\ref{galerkin}) 
with $mn=1430$ inner variables. 
We require a quadrature formula~(\ref{fcn-quadr}) to evaluate the nonlinear 
right-hand side approximately. 
We use a sparse grid quadrature, see~\cite{gerstner-griebel}, 
which is adapted from the Gauss-Legendre rule with level~3 
and $k=2441$ nodes. 
Initial values for the system~(\ref{galerkin}) are determined as in 
the previous example. 
Now the numerical solution with tolerances
$\varepsilon_{\rm rel} = 10^{-5}$, $\varepsilon_{\rm abs} = 10^{-6}$ 
yields 222~snapshots including the initial values in the 
time interval $[0,T]$.
On the one hand, the snapshots imply an approximation of the
coefficient functions in the truncated expansion~(\ref{pce-appr}) 
of the output voltage.
Figure~\ref{fig:transamp_coeff} illustrates the maximum coefficients 
occurring in the discrete time points.
We observe different orders of magnitudes, 
which indicates a potential for a sufficiently accurate low-dimensional 
approximation of the QoI.
On the other hand, a POD of the snapshots reveals the singular values 
within~(\ref{svd}) shown by Figure~\ref{fig:transamp_pod} (left). 
We obtain the associated projection matrices $T_{\rm l},T_{\rm r}$ 
for user-defined reduced dimensions.

Now we repeat the procedure for a stochastic collocation technique, 
where we apply the sparse grid quadrature from above. 
The auxiliary system~(\ref{collocation}) becomes a DAE with 
differential index one and $kn = 12205$ inner variables. 
Concerning the time integrations, the same tolerances are used 
as in the stochastic Galerkin approach.
The numerical solution of the initial value problem of~(\ref{collocation})
is interpolated onto 200 equidistant time points in $[0,T]$,
which yields 201~snapshots.
Alternatively, if the dynamical systems~(\ref{dae}) are solved separately
for the nodes of the quadrature (as done in Section~\ref{sec:scrapie}), 
then the interpolated snapshots cause failures within the transient simulation
of some reduced systems in this example.
Figure~\ref{fig:transamp_pod} (right) shows the dominating singular values 
from~(\ref{svd}) associated with the snapshots from the
collocation system~(\ref{collocation}).

Both the stochastic Galerkin system~(\ref{galerkin}) and the
collocation system~(\ref{collocation}) are reduced by the POD method now.
Initial value problems of FOMs and ROMs are solved with tolerances
$\varepsilon_{\rm rel} = 10^{-3}$, $\varepsilon_{\rm abs} = 10^{-6}$
in the following.
We consider the total time interval~$[0,3T]$, whereas the
snapshots are located in~$[0,T]$ only.
Sometimes the transient simulation of an ROM fails,
because the true dynamics is not captured. 
The reasons are that the dimension of an ROM is not large enough or
the snapshots do not reveal some required information.

\begin{figure}
\begin{center}
\includegraphics[width=6.5cm]{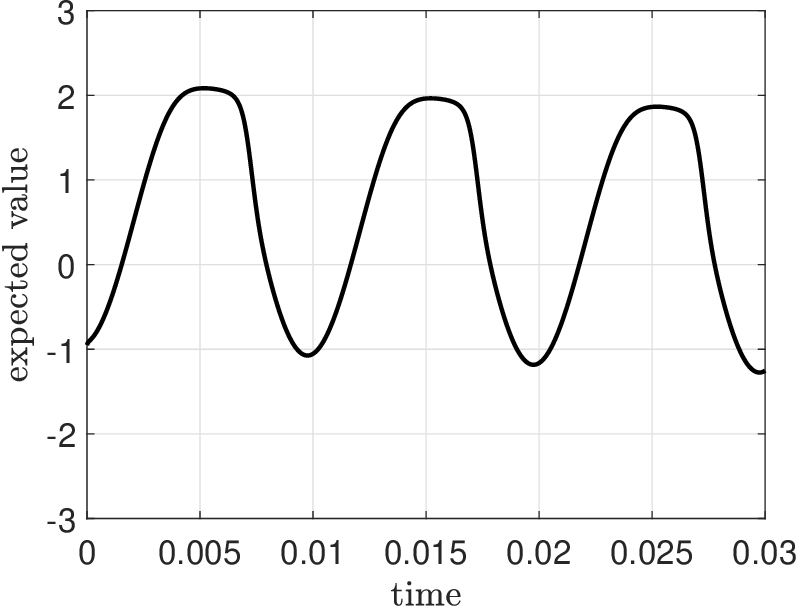}
\hspace{5mm}
\includegraphics[width=6.5cm]{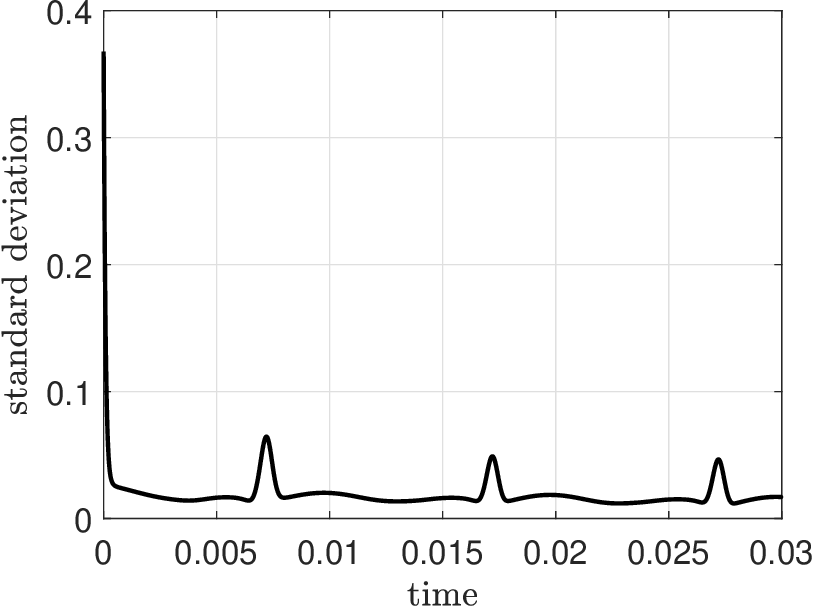}
\end{center}
\caption{Expected value (left) and standard deviation (right) for
  random output voltage computed by Monte-Carlo simulation.}
\label{fig:transamp_moments}
\end{figure}

A Monte-Carlo simulation generates a reference solution,
where initial value problems of the original DAEs~(\ref{dae})
are solved for $10^4$ samples of the random parameters.
The time integrations are done with high accuracy requirements
$\varepsilon_{\rm rel} = 10^{-6}$, $\varepsilon_{\rm abs} = 10^{-8}$.
Figure~\ref{fig:transamp_moments} depicts the computed expected value
as well as standard deviation.
We approximate the expected value as well as the variance of the QoI
using the FOMs and their ROMs.
These statistics are compared to the reference solution
in Table~\ref{tab:error}.
The maximum difference is determined on the time interval
$[0,T]$ of the snapshots and the longer time interval $[0,3T]$.
In the interval of the snapshots, the differences decay monotone
for increasing dimensions of the ROMs.
The ROMs do not improve from $r=25$ to $r=50$,
because the total error is already dominated by the quality of the snapshots
and the error of the time integration. 
In the longer interval, the monotonicity is not given for
increasing dimensions,
because the snapshots do not reproduce all required information.

\begin{table} 
  \caption{Maximum differences to reference solution
    in time intervals $[0,T]$ and $[0,3T]$ for expected value
    as well as variance of random output voltage
    in transistor amplifier example.}
\label{tab:error}
\begin{center}
  \begin{tabular}{lccccc}
    & & \multicolumn{2}{c}{expected value} & \multicolumn{2}{c}{variance} \\ 
    & ROM-dimension & $[0,T]$ & $[0,3T]$ & $[0,T]$ & $[0,3T]$ \\ \hline\hline
    Galerkin & FOM & $4.1 \cdot 10^{-3}$ & $4.1 \cdot 10^{-3}$ & $4.0 \cdot 10^{-4}$ & $4.0 \cdot 10^{-4}$ \\
             & $25$ & $4.5 \cdot 10^{-2}$ & $4.5 \cdot 10^{-2}$ & $8.3 \cdot 10^{-3}$ &  $8.3 \cdot 10^{-3}$ \\
             & $50$ & $4.5 \cdot 10^{-2}$ & $4.5 \cdot 10^{-2}$ & $8.3 \cdot 10^{-3}$ & $8.3 \cdot 10^{-3}$ \\ \hline
    Collocation & FOM & $5.4 \cdot 10^{-3}$ & $5.4 \cdot 10^{-3}$ & $3.5 \cdot 10^{-4}$ & $3.5 \cdot 10^{-4}$ \\
             & $10$ & $7.3 \cdot 10^{-3}$ & $4.3 \cdot 10^{-2}$ & $1.1 \cdot 10^{-2}$ & $5.6 \cdot 10^{-2}$ \\
             & $25$ & $4.2 \cdot 10^{-3}$ & $4.2 \cdot 10^{-3}$ & $3.4 \cdot 10^{-4}$ & $3.4 \cdot 10^{-4}$ \\
             & $50$ & $4.2 \cdot 10^{-3}$ & $5.6 \cdot 10^{-3}$ & $3.3 \cdot 10^{-4}$ & $1.9 \cdot 10^{-3}$ \\ 
\end{tabular}
\end{center}
\end{table}

\begin{figure}
\begin{center}
Galerkin method \hspace{4cm} collocation method

\vspace{3mm}

\includegraphics[width=7cm]{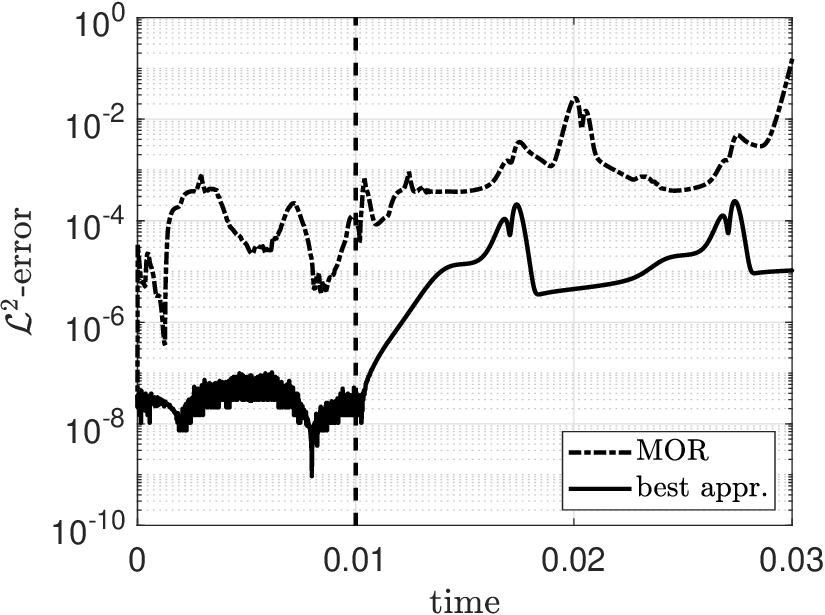}
\hspace{1mm}
\includegraphics[width=7cm]{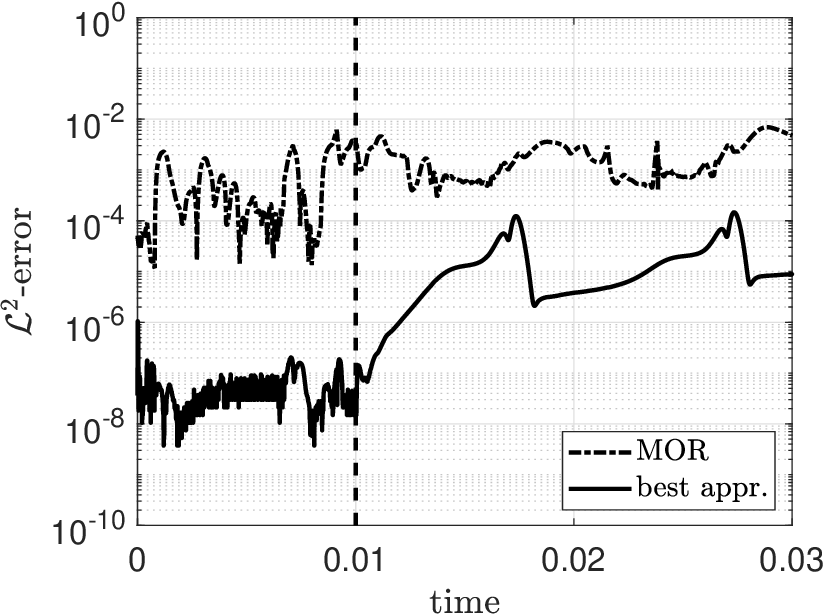}
\end{center}
\caption{$\ltwo$-errors for low-dimensional 
representations of QoI in the Galerkin method (left) 
and in the collocation technique (right) 
for reduced dimension~$r=25$ in transistor amplifier example.}
\label{fig:transamp_errors}
\end{figure}

We discuss QoIs from the reduced systems of dimension $r=25$ for both
Galerkin approach and collocation technique in the 
time interval $[0,3T]$ further.
We investigate the difference of an approximation~(\ref{qoi-appr}) 
from an ROM and a best approximation~(\ref{bestappr}) to a
solution~(\ref{qoi-method}) of an FOM.
Figure~\ref{fig:transamp_errors} depicts the computed
$\ltwo$-errors~(\ref{ltwo-error}) depending on time.
The behaviour of the error coincides in all four variants. 
The errors are relatively small in the first cycle, where
the snapshots are located. 
Yet the error increases at later cycles.
Of course the error of the best approximation is always smaller than the
MOR error.
In addition, the magnitudes of the $\ltwo$-errors agree for
Galerkin method and collocation method.
Hence the collocation technique becomes favourable in this example,
because the computational effort of a time integration is lower.

We note that computation work is not decreased in the used MOR for
this test example, because hyper-reduction is not included and thus
the nonlinear functions of the FOMs still have to be evaluated in
the ROMs. 
There is a potential for saving computing time by the usage of
hyper-reduction as mentioned in Section~\ref{sec:projection-mor}.



\section{Conclusions}
Orthogonal expansions were applied to the solution of 
random nonlinear dynamical systems. 
An MOR for the stochastic Galerkin system or an auxiliary system of 
a stochastic collocation method implied low-dimensional approximations 
of the expansions. 
On the one hand, a transient MOR method yields a low-dimensional
representation directly. 
On the other hand, the projection of a solution of the full-order model 
onto the subspace, 
which is identified by the MOR, generates a best approximation. 
Numerical simulations demonstrated that both the Galerkin method 
and collocation techniques are feasible to determine adequate
low-dimensional approximations. 
 However, a transient simulation of a reduced-order model
may become critical or fail
in the case of complex nonlinear dynamical systems. 
In contrast, the best approximation is robust and identifies 
an accurate low-dimensional approximation, while still 
the full-order model has to be solved.
 The strategy of hyper-reduction is required if the
proposed methods shall save computing time in the transient simulations.



\end{document}